\documentclass[a4paper,12pt,twoside]{article}
\usepackage{amssymb}
\usepackage[final]{lpretex}
\usepackage{title}
\usepackage{a4}
\usepackage{amscd}

\EndOfDump

\def\DHpartformat#1{(#1) }
\def\DHrefpart#1{(\DHRefpart{#1})}
\LBsetstyleof{partno}{arabic}

\DocVersion[Changed to LaTeX]{1}

\let\define\def

\def\A {{\mathbb A}}  \def\C {{\mathbb C}}
  \def\F {{\mathbb F}}

\def\GG {{\mathbb G}}   
  
\def\N {{\mathbb N}}  \def\P {{\mathbb P}} 
\def\Q {{\mathbb Q}} \def\R {{\mathbb R}}
 
\def\V {{\mathbb V}}  \def\X {{\mathbb X}}
\def\Z {{\mathbb Z}} 

\define \n {\mathbb N}
\define \z {\mathbb Z}
\define \q {\mathbb Q}
\define \PP {\mathbb P}

\def\sA {{\Cal A}}  
\def\sD {{\Cal D}} \def\sE {{\Cal E}} \def\sF {{\Cal F}}
  \def\sI {{\Cal I}}
  \def\sL {{\Cal L}}
 \def\sN {{\Cal N}} \def\sO {{\Cal O}}
 
  \def\sU {{\Cal U}}
  \def\sX {{\Cal X}}

\define \cN {\Cal N}
\define \cf {\Cal F}
\define \cg {\Cal G}
\define \cE {\Cal E}
\define \ce {\Cal E}
\define \cc {\Cal C}
\define \cV {\Cal V}
\define \cA {\Cal A}
\define \cK {\Cal K}
\define \cO {\Cal O}
\define \cF {\Cal F}
\define \cn {\Cal N}
\define \cI {\Cal I}
\define \sP {\Cal P}
\def\tA {\widetilde{\Cal A}}
\def\bA {\overline{\Cal A}}

\define \x {\xi}
\define \y {\eta}
\define \G {\Gamma}
\define \r {\rho}
\define \w {\omega}

\def \tU {\widetilde U}
\def \tZ {\widetilde Z}
\def\tX {\widetilde X}

\def \trho {\widetilde {\rho}}

\def \tpi {\widetilde{\pi}}

\def \tp {\widetilde{\mathbb P}}
\define \tH {\widetilde H}
\define \tG {\widetilde{\Gamma}}
\define \tW {\widetilde W}
\define \tF {\widetilde F}
\define \tm {\widetilde m}
\define \St {\widetilde S}
\define \Xt {\widetilde X}
\define \tS {\widetilde S}
\define \tpsi {\widetilde \psi}
\define \tL {\widetilde L}
\define \tE {\widetilde E}
\define \tl {\widetilde l}
\define \tA {\widetilde A}
\define \tom {\widetilde\omega}
\define \tT {\widetilde T}
\define \tB {\widetilde B}
\define \tf {\widetilde f}
\define \tsA {\widetilde{\sA}}

\define \tM {\widetilde M}
\define \tphi {\widetilde{\phi}}
\define \trho {\widetilde{\rho}}
\define \tR {\widetilde R}
\define \tp {\widetilde p}
\define \tq {\widetilde q}
\define \tc {\widetilde c}
\define \tsF {\widetilde {\sF}}
\define \tsN {\widetilde {\sN}}
\define \tsU {\widetilde {\sU}}
\define \th {\widetilde h}

\def\pd {\partial}

\def \Dx1 {\frac{\pd}{{\pd} x_1}}
\def \Dy1 {\frac{\pd}{{\pd} y_1}}
\def \Dz1 {\frac{\pd}{{\pd} z_1}}
\def \Dx2 {\frac{\pd}{{\pd} x_2}}
\def \Dy2 {\frac{\pd}{{\pd} y_2}}
\def \Dz2 {\frac{\pd}{{\pd} z_2}}

\def\q {\quad} 

\def\mapdiagr#1{\Big\searrow\rlap{$\raise 5pt\vbox{{\hbox{$\mkern -15mu\scriptstyle#1$}}}$}}   

\def\mapdiagl#1{\llap{$\raise 5pt\vbox{{\hbox{$\scriptstyle#1\mkern
-15mu$}}}$}\Big\swarrow}              

\def\Mapdiagr#1{\nearrow\rlap{$\lower 5pt\vbox{{\hbox{$\mkern
-15mu\scriptstyle#1$}}}$}} 

\def\Mapdiagl#1{\llap{$\lower 5pt\vbox{{\hbox{$\scriptstyle#1\mkern
-15mu$}}}$}\searrow} 

\def\Mapswr#1{\swarrow\rlap{$\lower 5pt\vbox{{\hbox{$\mkern
-15mu\scriptstyle#1$}}}$}}              

\def\Mapnwl#1{\nwarrow\rlap{$\lower 5pt\vbox{{\hbox{$\mkern
-15mu\scriptstyle#1$}}}$}}

\def \inj {\hookrightarrow}
\def \onto {\twoheadrightarrow}  

\define \Rhook {\hookrightarrow}

\def \half {\raise1pt\hbox{$\scriptstyle
        \frac{1}{2}\displaystyle$}}

\def \x{{\sl X}\llap{$\mkern -2mu {\scriptstyle -}$}}
\let\Spec\Sp

\def \Symm {\operatorname{Sym}}

\def \Pic {\operatorname{Pic}}

\define \Kod {\operatorname{Kod}}
\define \dimension {\operatorname{dim}}
\define \codim {\operatorname{codim}}
\define \contr {\operatorname{contr}}
\define \rk {\operatorname{rank}}
\define \im {\operatorname{im}}
\define \Mor {\operatorname{Mor}}
\define \Cl {\operatorname{Cl}}
\define \Hilb {\operatorname{Hilb}}
\define \degree {\operatorname{deg}}
\define \mult {\operatorname{mult}}
\define \Aut {\operatorname{Aut}}
\define \NS {\operatorname{NS}}
\define \Gal {\operatorname{Gal}}
\define \ch {\operatorname{char}}
\define \Jac {\operatorname{Jac}}
\define \Km {\operatorname{Km}}
\define \Sec {\operatorname{Sec}}
\define \Stab {\operatorname{Stab}}
\define \Br {\operatorname{Br}}
\define \inv {\operatorname{inv}}
\define \tr {\operatorname{tr}}
\define \Frob {\operatorname{Frob}}
\define \Symn {\operatorname{Sym}^n}
\define \Ev {\sE^\vee}
\define \ordp {\operatorname{ord}_p}
\define \Supp {\operatorname{Supp}}
\define \Ann {\operatorname{Ann}}
\define \disc {\operatorname{disc}}
\define \Lie {\operatorname{Lie}}
\define \embdim {\operatorname{embdim}}

\def\cone{\overline{NE}}

\def\onehalf{\frac{1}{2}}

\def\Chow{\operatorname{Chow}}

\def\NE{\operatorname{NE}}
\def\min{\operatorname{min}}
\def\Alb{\operatorname{Alb}}

\def\hod#1#2#3#4{\ensuremath{\if#30 H^{#2}({#1},{\cal O}_{#1}) \else 
 H^{#2}(#1,\Omega^{#3}\if\relax{#4}\relax_{#1}\else _{#1/#4}\fi)\fi}}

\begin{document}
\title{Perfect forms and the moduli space of abelian varieties}

\author{N. I. Shepherd-Barron}
\address{D.P.M.M.S.\\
 Centre for Mathematical Sciences\\
 Cambridge CB3 0WB\\
 U.K.}
\email{nisb@dpmms.cam.ac.uk}

\maketitle

Toroidal compactifications of
the moduli space $A_g$, or the stack $\sA_g$, of principally polarized
abelian $g$-folds have been constructed over $\C$ in \cite{AMRT} and over
any base in \cite{FC}. Roughly speaking, each such compactification corresponds 
to choosing a way of decomposing
the cone of real positive quadratic forms in
$g$ variables. The choice made here is the \emph{perfect cone} decomposition,
also called the \emph{first Voronoi} decomposition, which leads to the
compactification $A_g^{F}$. This carries divisor classes $M$, the
line bundle of weight $1$ modular forms, and $D$, the reduced boundary
$A_g^{F}\setminus A_g$. 
(Sometimes, but not in this paper, $M$ is denoted instead by $\omega$.)
It is easy to see that
$D$ is geometrically irreducible, and it follows (\cite{M} or \cite{F2})
that the classes $M$ and $D$ generate $\NS(A_g^{F})\otimes\Q$.
Here is the first main result of this paper, where a divisor class $E$ 
on a projective variety $X$ is \emph{nef}
if $E.C\ge 0$ for all curves $C$ on $X$:

\begin{theorem}\label{main1} 
$aM-D$ is nef if and only if $a\ge 12$ and ample if and only if $a>12$.
\noproof
\end{theorem}

This extends the the picture when $g=1$,
where there is a unique cusp form 
of weight $12$ and level $1$,
the discriminant, and it has no zeroes away from the cusp.
A better generalization would be the statement that
$12M-D$ eventually has no base points (that is, the complete 
linear system $\vert m(12M-D)\vert$ has
no base points for $m>>0$), but I can't prove that
except when $g\le 11$ and the base is $\C$.

One consequence is that the cone of curves
$\cone(A_g^{F})$ is the closed convex cone generated by curves
$C_1,C_2$, where $C_1$ is the closure of the set of points
$B\times E$, where $B$ is a fixed principally polarized abelian $(g-1)$-fold and 
$E$ is a variable
elliptic curve, and $C_2$ is any exceptional curve of the contraction
$A_g^{F}\to A_g^{Sat}$. Note that the properties of $\Delta$
just described
appear in this context as the formula $(12M-D).C_1=0$.
Another consequence 
is that for any value of $a> 12$,
the graded ring $\bigoplus_{n\ge 0}H^0(A_{g,\Z}^{F},\sO(n(aM-D)))$
of Siegel modular forms of weight $an$, vanishing along the principal cusp
to order at least $n$ (that is, of slope $a$) and with Fourier coefficients in $\Z$ 
is finitely generated over $\Z$.
 
The second main result derives from \ref{main1}. It concerns the nature of $A_g$,
now over $\C$, in the context of arbitrary quasi-projective varieties.

If $X$ is a quasi-projective complex algebraic variety of general type, then its 
\emph{canonical
model} (in the sense of Mori and Reid, not of Shimura) $X_{can}$ is a normal complex 
projective variety, birationally equivalent to $X$, with canonical singularities
and ample canonical class. It is known that $X_{can}$ exists if 
the dimension of $X$ is at most $3$, since this is when flips
are known to exist and terminate. However, even when the
canonical model is known to exist it is not always easy to find or describe
explicitly,
even for such concrete examples as Hilbert modular surfaces \cite{vdG}. 
The existence of $X_{can}$ is equivalent to the
canonical ring $\oplus H^0(\tX,\sO(nK_{\tX}))$ of $X$ 
being of finite type for some, or any,
smooth projective model $\tX$ of $X$. 
In any case, this ring is determined uniquely by the function 
field $\C(X)$ of $X$. There is also a local version of this notion: if $X$
is a normal complex variety, then its \emph{relative canonical model} is a
proper birational morphism $f:Y\to X$ such that $Y$ has canonical singularities
and $K_Y$ is ample relative to $f$. Again, this is uniquely determined by $X$
if it exists.

Freitag \cite{F1}, Tai \cite{T} and Mumford \cite{M}
have shown that, over $\C$, the coarse moduli space $A_g$ of 
principally polarized abelian $g$-folds, is of general type when $g\ge 7$.

\begin{theorem}\label{main2}\part[i]
$A_g^{F}$ is the relative canonical model
of the Satake-Baily-Borel compactification $A_g^{Sat}$ when $g\ge 5$.

\part[ii] When $g\ge 12$, the canonical model of $A_g$ exists and equals $A_g^{F}$.

\part[iii] $A_{11,can}$ exists and arises as the contraction of a certain extremal
ray in the cone of curves $\cone(A_{11}^F)$.
\noproof
\end{theorem}

Deriving \ref{main2} from \ref{main1} requires
the local statement that $A_g^{F}$ has canonical 
singularities over $\C$ when $g\ge 5$. For $A_g$ this is due to Tai \cite{T};
there are hints given there that this is true for $A_g^{F}$, but it seemed to be 
worth making it explicit.
From the formula $K_{A_g^{F}}=(g+1)M-D$,
the ampleness of $K_{A_g^{F}}$ for $g\ge 12$ is immediate. 

Namikawa \cite{N} has already raised the question of
the geometric meaning of the first Voronoi compactification
(which is also known as the \emph{perfect cone compactification}).
Alexeev \cite{A} has shown that, over any base, the \emph{second} Voronoi 
compactification $A_g^{Second}$ has an interpretation as a moduli space of generalized 
principally polarized abelian varieties; more precisely,
he defines and solves a moduli problem for generalized ppav's, shows that the moduli
space is proper and then shows that the second Voronoi compactification is an irreducible 
component of this space. By contrast, there is no reason to believe that $A_g^{F}$
is a moduli space; it is not even clear that there is an equidimensional 
family of projective 
schemes over it that extends the universal family of abelian varieties.
Another thing is that the natural birational equivalence between 
$A_g^{F}$ and $A_g^{Second}$
is regular in neither direction, because \cite{ER} when $g\ge 6$
neither of the two Voronoi decompositions is a refinement of the other.

Hulek and Sankaran \cite{HS02} raise the issue of describing the nef cone of 
$A_g^{Second}$. This seems much harder.
In \cite{HS04} they have already
described the nef cones of $A_g^{F}$ and $A_g^{Second}$ for $g\le 4$ and, 
which is most relevant to this paper, the nef cone of Mumford's
partial compactification $A_g^{part}$, the open subvariety lying over the open
subvariety $A_g\coprod A_{g-1}$ of $A_g^{Sat}$. That is, they deal
with complete curves in $A_g^{part}$. (This inverse image
is the same for all toroidal compactifications.)
\medskip 

I am grateful to Tom Fisher for his help in elucidating one of Tai's calculations
and to Klaus Hulek for his valuable remarks.
\bigskip
\begin{section}{Curves on the first Voronoi compactification.}\label{voronoi}
\medskip
Fix a copy $X_g$ of $\Z^g$. Denote by $X_g^\vee$ its dual and
$B(X_g)$ the lattice of symmetric bilinear
$\Z$-valued forms on $X_g$. Let ${\overline C}(X_g)$ denote the cone
in $B(X_g)\otimes\R$
of positive semi-definite forms with rational radical.
According to \cite{FC}, for every ring $R$, any \emph{basic} $GL(X_g)$-admissible
decomposition of 
${\overline C}(X_g)$ defines a smooth Deligne--Mumford stack $\tsA_{g,R}$
over $\Spec R$ that is a toroidal compactification 
of $\sA_{g,R}$, such that $\tsA_{g,R}=\tsA_g\otimes_{\Z}R$.
It also defines
a toroidal compactification
$\tA_{g,R}$ of the coarse moduli space $A_{g,R}$
over $\Spec R$. (When $R=\Z$ it is omitted.) 
The Satake compactification $A_{g,R}^{Sat}$ of $A_{g,R}$
is constructed in \cite{FC} as a blowing down of $\tA_{g,R}$.
However, as Alexeev 
remarks \cite{A}, it follows from this and the general principles of 
torus embeddings that any $GL(X_g)$-admissible
decomposition of ${\overline C}(X_g)$ defines a toroidal compactification over $\Spec R$.
(Alternatively, note that the only use made in \cite{FC} of the assumption that
the decomposition be basic is to ensure that the sheaf of $1$-forms with
logarithmic poles along the boundary is locally free; they use this to analyze
the Kodaira-Spencer maps of the semi-abelian schemes that they construct. Since this
sheaf is locally free 
for any toroidal scheme over $\Spec R$, the arguments of \cite{FC} carry over.)
Moreover, when $n$ is an integer such that $R$ contains $1/n$ and a primitive $n$th
root of unity, there are corresponding compactifications of the moduli objects 
$\sA_{g,n}$ and $A_{g,n}$ for ppav's with full level $n$ structure.
 
From the viewpoint of constructing a toroidal resolution of$A_{g,R}^{Sat}$, 
the choice of resolution over
the $0$-dimensional cusp determines what happens over all cusps,
in that a choice of admissible decomposition of ${\overline C}(X_g)$
determines an admissible decomposition of ${\overline C}(X_r)$ for
every quotient $X_g\onto X_r$.
There is a unique cusp (a copy of $A_{g-1}^{Sat}$) in $A_g^{Sat}$ of maximal dimension
and over this cusp there is a unique exceptional divisor $D$ in $\tA$,
which is generically the universal Kummer variety of dimension $g-1$.
This divisor corresponds to the primitive rank $1$ forms in ${\overline C}(X_g)$, 
which are all equivalent under $GL(X_g)$.

\begin{theorem}\label{decomposition} 
\part [i] \label{(2)} Taking the cones over the faces of 
the convex hull
of $(B(X_g)-\{0\})\cap {\overline C}(X_g)$ 
provides a $GL(X_g)$-admissible decomposition
of ${\overline C}(X_g)$.

\part [ii] \label{(1)}\label{barnes}
The convex hull of $(B(X_g)-\{0\})\cap {\overline C}(X_g)$
equals the convex hull of the set of primitive rank 
$1$ forms in ${\overline C}(X_g)$.
\begin{proof}
For \DHrefpart{i} see pp. 144-150 of \cite{AMRT}. In their notation, 
this decomposition is provided
by the \emph{perfect} co-core. Part \DHrefpart{ii} is the main result 
of \cite{BC}; in fact, they prove that every positive semi-definite form 
of rank at least two is in the
interior of the convex hull of the set of primitive rank one forms.
\end{proof}
\end{theorem}

Following Namikawa \cite{N} we shall refer to this decomposition of ${\overline C}(X_g)$ 
as the \emph{first Voronoi decomposition} or the \emph{perfect cone decomposition}
and the resulting compactification of $A_g$, or of any finite level cover
$A_{g,n}$ of it or of the stack $\sA_g$, as the first Voronoi compactification
or the perfect cone compactification,
denoted by $A_g^{F}$ or $\sA_g^{F}$. 
The reason is that 
$B(X_g)^\vee$ is naturally isomorphic to the lattice 
$Q(X_g)=\Symm^2(X_g)$ of
quadratic forms on $X_g^\vee$, so that the maximal faces of this decomposition
are level sets of some element $q$ of $Q(X_g)$ that is unique up to a scalar.
By definition, and the Barnes-Cohn result, these forms $q$ 
are exactly the \emph{perfect} ones.

By the construction given in
\cite{AMRT} or \cite{FC}, over a $0$-dimensional cusp, $\sA_g^{F}$
is formally isomorphic to the quotient stack $[{\overline E}/GL(X_g)]$, where 
${\overline E}$ is the locally finite torus embedding
$E\inj {\overline E}$ defined by taking $\mathbb X_*(E)=B(X_g)$ and then 
choosing the first Voronoi decomposition of ${\overline C}(X_g)$. 
Over an arbitrary cusp, 
$\sA_g^{F}$
is formally isomorphic to 
$[{\overline\Xi}/GL(X_{g-h})]$, where $X_{g-h}$ is a 
rank $g-h$ quotient lattice of $X_g$,
$E_1\inj {\overline E}_1$ is the torus embedding
described by replacing $X_g$ by $X_{g-h}$ in the description of $E\inj {\overline E}$
and ${\overline\Xi}$ is the contracted fibre product 
${\overline E}_1\times^{E_1}\Xi$, where
$\Xi$ is a certain
$E_1$-torsor over the $(g-h)$-fold fibre product of the universal
principally polarized abelian scheme $U_h\to\sA_h$.
We shall recall further details of this construction later on.

\begin{theorem}\label{first Voronoi} Assume that $g\ge 2$.

\part[i] The contraction $\pi:A_g^{F}\to A_g^{Sat}$ has a unique exceptional
divisor $D$. For every prime $p$, $D\otimes_\Z\F_p$ is absolutely irreducible.

\part[ii]\label{relample} 
$-D$ is ample relative to $\pi$.
 
\begin{proof} The exceptional divisors correspond to the $GL(X_g)$-orbits
of the vectors that form the
$1$-skeleton of the decomposition. The vectors in the $1$-skeleton
are the primitive rank $1$ forms, which are $GL(X_g)$-equivalent.

The projectivity criterion 
of \cite{FC}, V.5, shows that the inverse image of $-D$
on $A_{g,n}^1$ is ample relative to the contraction $A_{g,n}^1\to A_{g,n}^{Sat}$
for any $n\ge 3$. The result follows for $A_g^{F}\to A_g^{Sat}$ by taking quotients.
\end{proof}
\end{theorem}

The passage from a Deligne--Mumford
stack to its geometric quotient \cite{KM} does not always
commute with base change. The next lemma spells this out.

\begin{lemma}\label{specialization}
For any toroidal compactification $\tA_g$ as above and any field $k$, there is
a natural morphism $\tA_{g,k}\to\tA_{\Z}\otimes k$ that is birational and finite.

\begin{proof} This is an immediate consequence of the fact that for any $g$ and
over any algebraically closed field there are ppav's whose automorphism groupscheme
is precisely $\Z/2\Z$ (for example, take the Jacobian of a curve with trivial automorphism
groupscheme). So the natural morphism $\tA_{g,k}\to\tA_{\Z}\otimes k$,
which certainly exists, is an isomorphism over the locus where the automorphism groupscheme
is $\Z/2\Z$.
\end{proof}
\end{lemma}

\begin{theorem}\label{rho = 2} 
For every field $k$ the variety $A_g^{F}\otimes k$ is 
$\Q$-factorial and the $\Q$-vector spaces $\Pic(A_g^{F}\otimes k))\otimes\Q$
and $\Pic(\sA_g^{F}\otimes k))\otimes\Q$ 
are $2$-dimensional and generated by the classes $M$ and $D$.

\begin{proof} It follows from Corollary 1.6 of \cite{M} that the natural
homomorphism $H^2(A_{g,\C}^{F},\Q)\to\Pic(A_{g,\C}^{F})\otimes\Q$ 
is an isomorphism,
and the result then holds for any field of characteristic zero.

In general, choose a prime $l$ different from $\ch k$. We can assume $k$
to be algebraically closed. Then the specialization map
$$H^2(A_g^{F}\otimes{\overline\Q},\Q_l(1))\to H^2(A_g^{F}\otimes k, \Q_l(1))$$
is surjective. By \ref{specialization}, the natural homomorphism
$$H^2(A_g^{F}\otimes k, \Q_l(1))\to H^2(A_{g,k}^{F}, \Q_l(1))$$
is an isomorphism. Since
$\rho(A_{g,k}^{F})\ge 2$, from the existence of a contractible divisor,
we are done.
\end{proof}  
\end{theorem}

From now on we shall be careless in distinguishing between $A_g^{F}\otimes k$ and its
normalization $A^{F}_{g,k}$. The justification for this is \ref{specialization}.

Now suppose that the base is a field $k$.
Pick any curve $C_2\subset D$ that is contracted to a point in $A_{g}^{Sat}$.
Fix a principally polarized abelian $(g-1)$-fold $B$ and let $C_1$
be the closure in $A_g^F$ of the curve $\{B\times E\}$, 
where $E$ is a varying elliptic curve
and $B\times E$ is given the product principal polarization. The rest of this section
is devoted to showing that these curves generate the cone $\cone(A_g^{F})$. 

\begin{lemma}\label{classical}
$M.C_1=1/12$, $M.C_2=0$, $D.C_1=1$ and $D.C_2<0$.
\begin{proof} The formulae $M.C_1=1/12$ and $D.C_1=1$ follow from the existence
of the discriminant in the case $g=1$. The other formulae are obvious.
\end{proof}
\end{lemma}

\begin{lemma}\label{faces}
Suppose that $f:X\to Y$ is a non-constant morphism of normal projective varieties
over a field and that $C$ is a reducible
curve in $X$ that is contracted to finitely many points by $f$. Then 
for every component $C'$ of $C$ the
ray $\R_+[C']$ lies in the boundary $\partial\cone(X)$ of $\cone(X)$.
\begin{proof} Choose an ample divisor class $H$ on $Y$. Then
$\R_+[C']$ lies in $\cone(X)\cap (f^*H)^\perp$, which is contained in
$\partial\cone(X)$.
\end{proof}
\end{lemma}

\begin{lemma}\label{C_2 on boundary} $C_2$ lies in $\partial \cone(A_{g,k}^{F})$.

\begin{proof} This follows from \ref{faces} and
the choice of $C_2$ as a curve that is 
contracted by the morphism $A_{g,k}^{F}\to A_{g,k}^{Sat}$.
\end{proof}
\end{lemma}
\end{section}
\bigskip
\begin{section}{$12M-D$ is nef}\label{nef}
\medskip
In this section the base is a field $k$; there is no loss of generality
if this is taken to be algebraically closed.

To show that a $\Q$-divisor class $E$ on a projective variety is nef, 
it is enough to show
that $E$ is nef when restricted to the base locus of 
$\vert nE\vert$ for some $n$. Here is a sketch of how we do this 
for $12M-D$ on $A_g^{F}$.

\noindent (1) A consideration of Kummer varieties shows that 
the base locus of the linear system $\vert m(12M-D)\vert$
lies over the copy of $A_{g-2}^{Sat}$ in $A_g^{Sat}$.

\noindent (2) The irreducible components $Z$ of the part of $D$ lying
over cusps of higher codimension are closures of 
proper bundles $Z^0\to U_{g-r}^r$, where $U_{g-r}^r$ is
the $r$-fold fibre product of the universal abelian
scheme over $A_{g-r}$. The fibre of each $Z^0\to U_{g-r}^r$ 
is the boundary of some torus
embedding, and the action of this torus extends to $Z$. Then, via the Borel fixed 
point theorem, the torus action pushes the Chow point $[C]$ of a negative 
curve into the Chow point $[C_0]$ of a curve in some boundary stratum
that corresponds to some cone $\sigma$ that meets the interior
$C(X_r)$ of ${\overline C}(X_r)$, in such a way that $[C]-[C_0]$
is effective and is supported on curves $\Gamma$ that are contracted to points
in $A_g^{Sat}$. Since $M.\Gamma=0$ and $D.\Gamma<0$, by \ref{relample}, 
we can replace $C$ by $C_0$.
Moreover, from the nature of $\sigma$, this boundary stratum is 
a closure of a copy of $U^r_{g-r}$.

\noindent (3)  The normalized closure $\tU$ of this copy $U$
of $U^r_{g-r}$ is, after deleting a
closed subset of codimension 2 (here the argument depends on the fact
that we're dealing with the perfect compactification), a semi-abelian scheme 
$U_{g-r}^{r,part} \to A_{g-r}^{part}$. So when we pull $12M-D$ back
to $\tU$, it will follow that any negative curve in $\tU$ will lie in
$\tU\setminus U$ provided that $m(12M-D)$ when restricted to the open
$U$ has no base points.

\noindent (4) By considering various projections 
$U_{g-r}^{r,part}\to U_{g-r}^{part}$
that are determined by the geometry of the cone $\sigma$
(here again we rely on having the perfect compactification)
we reduce the problem of proving that the restriction of $m(12M-D)$ to $U$
has no base points to what was proved in (1) and conclude by induction.

Write ${D_g}$ or $\sD_g$ for the (irreducible) boundary divisor
in $A_g^{F}$ or $\sA_g^{F}$. 

\begin{theorem}\label{modify} 
Suppose that the base is an algebraically closed field $k$. 
For some $m>0$ the linear system $\vert m(12M-{D_g})\vert$
has no base points in $A_g^{part}$.

\begin{proof} 
First, assume that $\ch k$ is not $2$. 
We shall work at level $2$.

The inverse image $D_g^{part}$ in $A_g^{part}$ of $A_{g-1}$
is the locus of (torus) rank $1$ degenerations. (In fact, $A_g^{part}$
and $\sA_g^{part}$ are independent of the choice of toroidal compactification.)
The universal abelian scheme over $\sA_g$
extends to an equivariant relative compactification $\pi:\sU_g\to \sA_g^{part}$
of a semi-abelian scheme.
Each fibre $\pi^{-1}(x)$ over a geometric point $x$ of $\sD^0$ is constructed
from a $\P^1$-bundle over a principally polarized
abelian $(g-1)$-fold by identifying two disjoint sections
and is an equivariant compactification of the semi-abelian variety
that is the smooth locus.

Consider the full level $2$ version
$\sA_{g,2}^{F}$ 
of the stack $\sA_g^{F}$ and the corresponding open substack $\sA_{g,2}^{part}$.
Again the universal abelian scheme over $\sA_{g,2}$ extends to an equidimensional 
projective family $\pi_2:\sU_{g,2}\to \sA_{g,2}^{part}$. This time, 
a fibre $F$ over a geometric point of $\sD^{part}$ is the sum of two components
$F_1$ and $F_2$. Each $F_i$ is a copy of $X\times \P^1$, where $X$ is 
a principally polarized
abelian $(g-1)$-fold, $F_1\times\{0\}$ is identified with
$F_2\times\{\infty\}$ by a translation determined by a chosen point $P_0\in X$
and similarly for 
$F_1\times\{\infty\}$ and $F_2\times\{0\}$.

There is a symmetric line bundle $\sL$ on $\sU_{g,2}$ that
defines twice the principal polarization on the smooth fibres.
On a singular fibre $F$ it is of degree $1$ on each copy of $\P^1$ that
appears and defines twice the principal polarization on $X$.

At this point we need a lemma.

\begin{lemma}\label{no base points}
For any $x\in A_{g,2}^{part}$
the specialization map 
$$\phi_x:\pi_{2*}\sL \to H^0(\sU_x,\sL_x)$$
is surjective and $\sL_x$ is generated by its sections. 

\begin{proof} The surjectivity
of $\phi_x$ is a consequence of the vanishing of $H^i(\sU_x,\sL_x)$ for $i>0$
and the standard base change results. 

For smooth fibres the absence of base points is well known.
For any fibre the space $H^0(\sU_x,\sL_x)$ is a representation of
the Heisenberg group $Heis_{g,2}$ (the central extension
of $\Z/2\Z^g\times\mu_2^g$ by $\mu_2$ determined by the Weil pairing)
and for a singular fibre $F=\sU_x$ in which $X$ appears as a component of
the singular locus, the restriction map
$H^0(F,\sL_x)\to H^0(X,(\sL_x)\vert_X)$ is equivariant for the
subgroup $Heis_{g-1,2}$ of $Heis_{g,2}$. Since $H^0(X,(\sL_x)\vert_X)$
is an irreducible representation of $Heis_{g-1,2}$, 
restriction is surjective. (If it were zero, then it would be zero on the other 
component of the singular locus, which is impossible since $\sL_x$
has degree $1$ on each $\P^1$.)
In particular, there are no base points in the singular locus of $F$.

Suppose that $P$ is a base point in the smooth locus. Consider the
copy $\Gamma$ of $\P^1$ in $F$ that contains $P$. Then there is a copy of $\mu_2$
in $Heis_{g,2}$ that preserves $\Gamma$ and acts on the restriction
of $H^0(F,\sL_x)$ to $\Gamma$; it follows that $\Gamma$ lies in the
base locus, and the lemma is proved.
\end{proof}
\end{lemma}

Since $Heis_{g,2}$ rigidifies the projective space
$\P(H^0(\sL_x))=\P^N$ (which is independent of the geometric point $x$),
where $N=2^g-1$. Therefore there is a morphism
$Km:\sA_{g,2}^{part}\to \Chow(\P^N)$ that sends each point $x$ to the cycle
that is the image, counted with appropriate multiplicity, under
the $2\theta$ linear system, of the scheme $f_2^{-1}(x)$. 
If $x$
corresponds to an irreducible principally polarized abelian variety, 
then $Km(x)$ is the associated
Kummer variety, counted with multiplicity $1$ (\cite{LB}, Th. 8.1, 
p. 99; the assumption there that the case field is $\C$ is unnecessary
for their argument), so that $Km$ is radicial (that is, generically an
\'etale homeomorphism)
onto its image. However, $Km$ is constant along the inverse image in
$\sA_{g,2}^{part}$ of the curve $C_1$; this is the statement that 
the Kummer variety of an elliptic curve is $\P^1$, so that
for $x=[E\times B]\in C_1$ the image of $\Km(E\times B)$ is just
$\P^1\times\Km(B)$ embedded in $\P^N$ via the Segre
embedding of $\P^1\times\P^{2^{g-1}-1}$. 
Since $Km$ factors through the geometric quotient $A_{g,2}^{part}$
of $\sA_{g,2}^{part}$,
we have a non-constant morphism defined on $A_{g,2}^{part}$
that collapses the inverse image of $C_1$.
To get the morphism that we want on $A_g^{part}$, take
quotients by the finite group
$Sp_{2g}(\F_2)$: acts on $\P^N$, so on $\Chow(\P^N)$, and gives rise
to a commutative diagram
$$
\begin{CD}
\sA_{g,2}^{part} @>{H}>> \Chow(\P^N)\cr
 @VVV  @VVV\cr
\sA_g^{part} @>{h}>> \Chow(\P^N)/Sp_{2g}(\F_2).
\end{CD}
$$
Observe that the morphism $h$ factors through a morphism $H$
defined on $A_g^{part}$ that contracts $C_1$. 

Regard $H$ as a rational map on $A_g^{F}$. Since $\rho(A_g^{F})=2$
and $A_g^{F}$ is regular,
$H$ is defined by some linear system $\vert aM-bD\vert$.
Since $H$ contracts $C_1$ and $M.C_1=\frac{1}{12}$, $D.C_1=1$,
by \ref{classical},
$am-bD$ is proportional to $12M-D$.

Now suppose that $\ch k=2$.
The proof follows the same lines, except that we consider
a full level $3$ structure,
the linear system $\vert3\Theta\vert$ and the action of the level $3$ 
Heisenberg group on (the dual of) this projective space. To carry this out 
requires the introduction of a theta level structure, since otherwise there 
is no line bundle $\sO(3\Theta)$.

According to \cite{FC}, p. 132 \emph{et seq.}, there is, over any ring $R$
that contains $\zeta_3$ and in which $3$ is invertible,
a stack $\sN_{g,3}$, finite and flat over $\sA_{g,3}$, with a universal
abelian scheme $\sU^0\to\sN_{g,3}$ with a full level $3$ structure
and a symmetric ample line bundle $\sO(\Theta)$
that induces a principal polarization of $\sU_g^0$. 
Now suppose that $R$ is a field. Then, in the notation of \emph{loc. cit.}
(except that we write $L$ instead of $X$),
the perfect cone decomposition, for
a choice of $\rho\in L$, of 
the cone ${\overline C}(X_g)$, where now $Q_\rho(L)$ is the lattice
instead of $B(X_g)$, so that the perfect cone decomposition
is the convex hull of the rank $1$ forms in $Q_\rho(L)$ instead of $B(X_g)$, 
determines a normal toroidal compactification
$\sN^1_{g,3}$ of any component $\sN_{g,3}^{norm}$ of the normalization  of
$\sN_{g,3}$. Moreover, the abelian scheme $\sU^0\to \sN_{g,3}$
extends to a semi-abelian scheme over $\sN_{g,3}^{part}$ that possesses
an equivariant compactification $\sU\to\sN_{g,3}^{part}$ on which there
is a line bundle $\sL$ that is a relative $\sO(\Theta)$.

A slight modification of the proof of \ref{no base points} shows that
the line bundle $\sO(3\Theta)$ is relatively very ample over 
$\sN_{g,3}^{part}$.
When $X$ is, as a $g$-dimensional principally polarized abelian variety, 
of the form $X=E\times B$ with $E$ an elliptic curve,
the intersection of the quadrics containing the $\vert3\Theta\vert$-image of $X$
is $\P^2\times B'$, where $B'$ is the intersection of the
quadrics that contain the $\vert3\Theta\vert$-image of $B$. In particular,
the $\vert3\Theta\vert$-image of $X$
lies in the image $\Sigma$ of the Segre embedding 
of $\P^2\times\P^{3^{g-1}-1}$
in $\P^{3^g-1}$. On the other hand, if $X$ is irreducible and $g\ge 3$, then its
$\vert3\Theta\vert$-image
cannot lie in $\Sigma$, for then $X$ would possess a non-constant
morphism to $\P^2$. (The case where $g=2$ is left to the reader.)
Finally, instead of $\Chow(\P^N)$ we take
the Grassmannian of $M$-dimensional subspaces of the vector space $V$ of quadrics
in $\P^{3^g-1}$, where $M=\dim H^0(\P^{3^g-1},\sI_{X/\P^{3^g-1}}(2))$. 
(Note that $M$ is independent of the principally polarized abelian $g$-fold $X$,
since the $\vert3\Theta\vert$-image of $X$
is projectively normal
and $\dim H^0(X,\sO(6\Theta))=6^g$, so is independent of $X$.)
This space $V$ is rigidified
by the finite Heisenberg group 
determined by a full level $3$ structure, and now we argue as in
the first case.
\end{proof}
\end{theorem}

\begin{corollary} Any curve $C$ in $A_g^{F}$ with $(12M-{D_g}).C<0$ lies
in $A_g^{F}\setminus A_g^{part}$ and has $1$-dimensional image in $A_g^{Sat}$.
\begin{proof} The first part follows from \ref{modify} and the second
from the fact that $-{D_g}$ is ample relative to $A_g^{F}\to A_g^{Sat}$.
\end{proof}
\end{corollary}

So we must consider what happens higher up in the boundary. For this,
recall further details of the toroidal construction.

We have a lattice $X_g=\Z^g$ with quotient lattices
$X_g\onto X_{g-1}\onto\cdots\onto X_0=0$. 
Put ${\overline C}(X_r)={\overline C}_r$.
Then ${\overline C}_r$ is a subcone of ${\overline C}_{r+1}$.
Fix an admissible decomposition $\{\sigma\}$ of ${\overline C}(X_g)$; 
this determines an
admissible decomposition of ${\overline C}(Y)$ for every $X_g\onto Y$.
Over $A_{g-r}\subset A_g^{Sat}$, the picture is this (\cite{FC}, p. 105):
let $U_{g-r}\to \sA_{g-r}$ be universal, $U_{g-r}^r=Hom(X_r,A)$,
a copy of the $r$-fold fibre product of $U_{g-r}\to \sA_{g-r}$
and take a certain 
torsor $\Xi\to U_{g-r}^r$ under the $\frac{r(r+1)}{2}$-dimensional
torus $E=E_r$ with 
$\X_*(E_r) =B(X_r)$; the characterization of this torsor in terms of
the $\GG_m$-bundles over
$U_{g-r}$ associated to given characters of $E_r$
is given in the last paragraph of \emph{loc.cit.} and will be needed
later in \ref{product}.
Take the locally finite
torus embedding $E_r\inj{\overline E_r}$ associated to the admissible decomposition,
with boundary $\partial E_r={\overline E_r}\setminus E_r$. Notice
that the irreducible components of $\partial E_r$ correspond to the minimal
cones in the decomposition of ${\overline C}_r$ that meet the relative interior $C(X_r)$.
Then the inverse image of the locally closed subvariety $A_{g-r}$
of $A_g^{Sat}$ in the toroidal compactification $A_{g,\{\sigma\}}$
is the quotient by $GL(X_r)$ of the associated bundle
$\Xi\times^{E_r}\partial E_r\to U_{g-r}^r$. Note that this quotient by $GL(X_r)$ 
does exist in the category of schemes locally of finite type, 
and the analogous statement holds
at the level of stacks. From this, and consideration
of the embedding ${\overline C}(X_r)\inj{\overline C}(X_{r+1})$, the following lemma,
except maybe for \ref{closure}, whose significance is that the $Z$
that appears there is proper, is clear. 

\begin{lemma}\label{fixed point}
\part[i]\label{closure}
At level $n\ge 3$, the torus $E_r$ acts on each component $Z$ of the inverse image
of $A_{g-r,n}^{Sat}$ in $A_{g,n,\{\sigma\}}$ such that 
the inverse image of $A_{g-r,n}$ is finitely stratified with each stratum
being a bundle over $U_{g-r,n}^r$ whose fibre is a quotient of $E$.

\part[ii] The strata correspond to equivalence classes, under the principal
congruence subgroup of level $n$ in $GL(X_r)$, of the cones in $\{\sigma\}$ that 
lie in ${\overline C}(X_r)$ and meet $C(X_r)$.

\part[iii]\label{strata} The strata corresponding to maximal cones in
${\overline C}(X_r)$ are copies of $U_{g-r,n}^r$.

\begin{proof} For \DHrefpart{i}, take the morphism
$Z\to A_{g-r,n}^{Sat}$; we have seen that $E_r$ acts on the 
inverse image $Z^0$ of $A_{g-r,n}$. It is also clear, from considering
the embedding ${\overline C}(X_r)\inj{\overline C}(X_{r+s})$, that $E_r$
acts on the formal completion of $Z$ along its fibre over each 
copy of $A_{g-r-s,n}$ that is a cusp in $A_{g-r,n}^{Sat}$ 
compatibly with its action on $Z^0$. Since $Z$ is normal,
it is enough to show that if 
$Z$ is a normal $k$-variety
$U$ is open in $Z$, 
$F=Z\setminus U$ and $E$ is a $k$-torus
that acts compatibly on $U$ and 
on the $F$-adic completion $\widehat Z$ of $Z$, then $E$ acts on $Z$.

By assumption, there is a commutative diagram
$$
\begin{array}{ccc}
{E\times{\widehat Z}} & \to &{\widehat Z}\\
\downarrow &\ & \downarrow\\
E\times Z & -\to &  Z
\end{array}
$$
where the products are over $k$ and the horizontal arrows are the actions;
the lower one is a rational map whose base locus lies in $E\times F$.
This shows that the base locus of the rational map $E\times Z -\to Z$ disappears
after making the cover $E\times{\widehat Z}\to E\times Z$
of normal schemes; since this is cover is faithfully flat over $E\times F$, 
the base locus is empty.

The rest of the lemma has already been proved in the
discussion just preceding it.
\end{proof}
\end{lemma}

Note that the maximal cones referred to in \ref{strata} necessarily
meet $C(X_r)$.

Now suppose that $\{\sigma\}$ is the perfect cone decomposition of 
${\overline C}(X_g)$.
Its defining property is that for any cone $\sigma$ in this decomposition that
meets $C(X_r)$, there is a positive definite quadratic form in $r$ variables
the squares of whose minimal vectors span $\sigma$. This form is not unique,
even up to scalar multiplication, unless $\sigma$ is of maximal dimension
$\frac{r(r+1)}{2}$. In this case the form is perfect, essentially by definition
of that word. Conversely, if $q$
is a positive semi-definite form in $g$ variables, then the squares of its
minimal vectors span a cone in the decomposition; this repeats \ref{barnes}.


\begin{corollary}\label{negative}
If there is a curve $C$ in $A_g^{F}$ with $(12M-{D_g}).C<0$, then
for some $r\ge 2$ there is such a curve lying in the closure 
${\overline U}_{g-r}^r$ of 
the copy of $U_{g-r}^r$ that corresponds, under \ref{strata}, to 
some cone $\sigma$ in the perfect cone decomposition of ${\overline C}(X_g)$
that lies in ${\overline C}(X_r)$ and is maximal there.

\begin{proof} Work at level $n$, with $n\ge 3$ and prime to
$\ch k$. There is a natural cover $\pi:A_{g,n}^{F}\to A_g^F$
with $\pi^*{D_g}=nD_{g,n}$, 
where $D_{g,n}$ is the reduced sum of the boundary divisors
in $A_{g,n}^{F}$. So
having such a curve is equivalent, at level $n$,
to having a curve $C$ in $A^F_{g,n}$ with $(\frac{12}{n}M-D_{g,n}).C<0$.

We know that $C$ is in $A_{g,n}^{F}\setminus A_{g,n}^{part}$. So there is 
a unique value of $r$ such that the image of $C$ in $A_{g,n}^{Sat}$
lies in 
$A_{g-r,n}^{Sat}\setminus A_{g-r-1,n}^{Sat}$. Then the result
follows from \ref{fixed point} by applying the Borel fixed point theorem 
to the action of $E_r$ on the Chow scheme
of curves in the inverse image of $A_{g-r,n}^{Sat}$ in $A_{g,n}^{F}$,
as was explained at the start of this section.
\end{proof}
\end{corollary}

From now on, assume that there is a curve on which $12M-{D_g}$ is negative.
Take a cone $\sigma\subset B(X_r)\otimes\R$ and corresponding closure 
${\overline U}_{g-r}^r$
that contains such a negative curve, as
provided by \ref{negative}. Suppose that $q=q_\sigma\in Q(X_r)$
is the perfect form (unique up to a scalar) that defines $\sigma$. 
Let $\tU_{g-r}^r$ be the normalization of ${\overline U}_{g-r}^r$. Then,
in order to get further control on these negative curves,
it is enough to show that the complete linear system given by
some large multiple of the pullback of the divisor class $12M-{D_g}$ to
$\tU_{g-r}^r$ has no base points in the open subvariety $U_{g-r}^r$.

\begin{proposition}\label{good}
Suppose that $\tU_{g-r}^r$ is the normalization
of ${\overline U}_{g-r}^r$. Then the complement $\tU_{g-r}^r\setminus U_{g-r}^r$
is an irreducible divisor.

\begin{proof} This is a consequence, via the usual rules governing
the construction of torus embeddings, of the statement that 
the cones $\tau$ in $\overline C(X_{r+1})$ such that 
\smallskip

\noindent $(1)$
$\sigma$ is a proper face of $\tau$ and 

\noindent $(2)$ $\tau$ is minimal with
respect to this property 
\smallskip

\noindent form a single orbit under the subgroup
of $GL(X_{r+1})$ that preserves the quotient homomorphism
$X_{r+1}\onto X_r$ and acts trivially on $X_r$. 

To prove this, suppose that $f=f(x_1,\ldots,x_r)$ is a perfect form in $r$ 
variables and minimum value $a$. Then the 
set $\min(f_1)$ of minimal vectors of $f_1:=f+ax_{r+1}^2$
is given by 
$$\min(f_1)=\min(f)\cup\{\pm(0,\ldots,0,1)\},$$ 
so that, by the defining property
of the perfect cone decomposition, $f_1$ defines one of the
cones $\tau$ in question. 

To show that this construction accounts for all the $\tau$,
suppose that $\tau'$ is any one of the cones in question. There is a form
$F=F(x_1,\ldots,x_r,x_{r+1})$ the squares of whose minimal vectors span $\tau'$,
again by the defining property of the perfect cone decomposition.
Since $\sigma$ is spanned by a subset $S$ of these squares, the minimality
of $\tau'$ ensures that it is spanned by $S$ and just one more square.
So we have accounted for all the $\tau$.
\end{proof}
\end{proposition}

\begin{corollary}\label{codimension 2}
There is an open substack $\tsU_{g-r}^{r,part}$ of
$\tsU_{g-r}^r$ whose complement has codimension $2$ and that contains $\sU_{g-r}^r$ 
and is a semi-abelian scheme over $\sA_{g-r}^{part}$ extending
$\sU_{g-r}^r \to \sA_{g-r}$.
\begin{proof} This follows directly from the restatement of \ref{good}
in terms of stacks.
\end{proof}
\end{corollary}

Here is a summary of some of these objects and morphisms:
$$U_{g-r}^{r,part}\inj\tU_{g-r}^{r,part}\inj \tU_{g-r}^r\to {\overline U}_{g-r}^r
\inj D_{g}\inj A_{g}^{F},$$
where the embeddings are open except for $D_{g}\inj A_{g}^{F},$
which is closed.
Let $\alpha_{g,r}:U_{g-r}^{r,part}\to A_{g}^{F}$ denote the composite.
When $r=1$, ${\overline U}_{g-r}^r =D_{g}^{F}$. In addition,
there is a $0$-section $N_{g,r}^0 \subset U_{g-r}^r$.
Let $N_{g,r}$ denote the normalized closure of $N_{g,r}^0$
in $\tU_{g-r}^r$ and $N_{g,r}^{part}=N_{g,r}\cap \tU_{g-r}^{r,part}$.

\begin{lemma}\label{natural} $N_{g,r}$ is naturally isomorphic to $A_{g-r}^{F}$.
\begin{proof} First, the closure of $\sA_{g-r}\times \sA_r$ in $\sA_g^{F}$
is just $\sA_{g-r}^{F}\times \sA_r^{F}$; this is the statement that if $q_1$
is a positive quadratic form in $g-r$ variables, $q_2$ a positive quadratic
form in $r$ further variables and $q_1,q_2$ have equal minimum values,
then $\min(q_1+q_2)=\min(q_1)\cup \min(q_2)$.

Then inside $A_g^{F}$, $N_{g,r}$ is the limit of the locus
$A^{F}_{g-r}\times\{[B]\}$, where $[B]\in A_r$ and $[B]$ tends
to a point in $D_r$
in a way specified by the choice of $\sigma$ as a maximal cone in the perfect 
cone decomposition of ${\overline C}_r$. 
\end{proof}
\end{lemma}

The next thing is to understand the restriction of the class $[{D_g}]$ to 
its subvariety $N_{g,r}$ in two different ways.

\begin{lemma} The restriction of ${D_g}$ to $N_{g,r}$ is linearly
equivalent to $D_{g-r}$ via the identification $N_{g,r}=A_{g-r}^{F}$
of \ref{natural}.
\begin{proof} This follows at once from the limiting description
of $N_{g,r}$ given in the proof of \ref{natural}.
\end{proof}
\end{lemma}

For another view of this restriction, look again
at the embedding ${\overline U}_{g-r}^r\inj D_{g}$
determined by the cone $\sigma$.
Along this subvariety, $D_{g}$ is the union of branches 
$\delta_1,\ldots,\delta_p$ of which
${\overline U}_{g-r}^r$ is the intersection. 
At the stack level, these branches correspond, once a basis of $X$ has been chosen,
so that $X$ is identified with $X^\vee$, to 
the squares $x_1^2,\ldots,x_p^2$ of the
minimal vectors of the perfect form $q_\sigma$ in $r$ variables 
that determines $\sigma$,
or, equivalently, to the rank $1$ quadratic forms that span $\sigma$.
The number $p$ is half the kissing number of $q_\sigma$.

Denote by $L_{g-r}$ the symmetric divisor class on 
$\sU_{g-r}\to \sA_{g-r}$ that defines $(-2)$ times the principal polarization
on each fibre.
Let $\gamma:\sU_{g-r}^{r,part}\to \sA_{g-r}^{part}$
denote the structural morphism.

\begin{lemma}\label{H1}
Two divisor classes on $\sU_{g-r}^{r,part}$ or 
$U_{g-r}^{r,part}$ are equivalent, modulo torsion,
if they are equivalent on the zero-section and on the generic fibre.
\begin{proof} Since the complement of $\sU_{g-r}^r$ both in
$\sU_{g-r}^{r,part}$ and in the normalized closure $\tsU_{g-r}^r$ is an irreducible 
divisor, it is enough to prove the result for classes in $H^2$ on
$\sU_{g-r}^r$ over $\C$.

In this case it is enough, by the Hochschild--Serre spectral sequence,
to show that $H^1(Sp_{2g},M)$ is torsion, where $M=\Z^{2g}$ is the 
standard representation. As explained on p. 135 of \cite{BMS}, 
this follows from the linear reductivity
of the algebraic group $Sp_{2g}$ over $\Q$.
\end{proof}
\end{lemma}

\begin{proposition}\label{product} 
\part [i] The boundary divisor $D_{g-r+1}$ in $A_{g-r+1}^{F}$
cuts out $L_{g-r}$ via the identification
of $\sU_{g-r}$ with an open substack of $\sD_{g-r+1}$.

\part [ii] There are positive rational numbers
$c,b_i$ and projections
$\pi_i:U_{g-r}^{r,part}\to U_{g-r}^{part}$ over $A_{g-r}^{part}$,
for $i=1,\ldots,p$,
with composites 
$\rho_i=\alpha_{g-r+1,1}\circ\pi_i:U_{g-r}^{r,part}\to A_{g-r+1}^{F}$
such that $\sum b_i\rho_i^*(12M-D_{g-r+1})+c\gamma^*(12M-D_{g-r}^{part})$
is numerically equivalent to 
$\alpha_{g,r}^*(12M-D_{g})$.

\begin{proof} 
First, take the affine 
torus embedding $E\inj {\overline E}$ that corresponds to $\sigma$.
In particular, $B(X_r)=\X_*(E)$.

Since the $1$-skeleton of $\sigma$
is generated by vertices that lie in the same affine hyperplane,
there is a morphism $f:(E\inj {\overline E})\to (\GG_m\inj \A^1)$ of torus embeddings
such that the divisor $f^{-1}(0)$ is an integral multiple of the
reduced boundary divisor $\partial{\overline E}={\overline E}\setminus E$.
Fix a basis of $X$ and then identify $X$ with $X^\vee$. 
The lattice homomorphism
$\th:\Z^p=\sum \Z.e_i \to B(X_r)$, where $\th(e_i)=x_i^2$ and the $x_i$ 
are the minimal vectors of $q$, gives a morphism
$h:(\GG_m^p\inj \A^p)\to (E\inj {\overline E})$ of torus embeddings. Composing
this with $f$ shows that such that
$h^{-1}(\partial {\overline E})$ is the sum of the co-ordinate hyperplanes
in $\A^p$.

Now look at various associated bundles over $\sU_{g-r}^r$
constructed from the $E$-torsor $\Xi\to \sU_{g-r}^r$ described on p. 105
of \cite {FC}.
Also, consider $\sU_{g-r}^r$ as an open subscheme of $\tsU_{g-r}^r$,
so inside $\sA_g^F$. This copy $\sU_{g-r,origin}^r$ of $\sU_{g-r}^r$ appears inside
these associated bundles as the origin. Moreover, the associated $\A^p$-bundle 
$\V\to \sU_{g-r}^r$ that we get is the direct sum of certain line bundles 
$H_1,\ldots,H_p$. There is a morphism $h:\V\to {\overline \Xi}$
of bundles over $\sU_{g-r}^r$, and
the restriction of the $\Q$-divisor class $h^{-1}(\partial {\overline \Xi})$ to 
$\sU_{g-r,origin}^r$ 
is a rational multiple of $\sum c_1(H_i)$.

Note that, by the nature of the toroidal compactification $A_g^F$,
the $\Q$-divisor classes ${\overline \Xi}$ and $D_g$ have the same restriction
to $\sU_{g-r,origin}^r$. So the restriction of $D_g$ to $\sU_{g-r,origin}^r$
is the same as the restriction of $\sum\delta_i$ and is a rational
multiple of $\sum c_1(H_i)$, say $s\sum c_1(H_i)$.

We must describe these bundles $H_i$.
There is an identification
of $X_r\otimes (\sA^t)$ with $(X_r\otimes \sA)^t$, so that if 
$\lambda :\sA\to \sA^t$ is the given principal polarization, then 
$1_{X_r}\otimes\lambda: X_r\otimes \sA\to X_r\otimes (\sA^t)$ 
is a principal polarization.
It also identifies $X_r$ with $X_r^\vee$ and embeds 
$B(X_r)$ into its dual $Q(X_r)=\Symm^2(X_r)$.
Moreover, the elements $\th(e_i)$ of $B(X_r)$ are then characters of $E$
and the $\GG_m$-bundle 
$\Xi\times^{E,e_i}\GG_m\to \sU_{g-r}^r$
is $H_i\setminus\{0\}$.
Since $\th(e_i)$ is, as an element of $Q(X_r^\vee)$, the square of a primitive
vector in $X_r$, it follows from the description given on p. 105 of \cite{FC}
that $H_i$ is the pullback of the inverse Poincar\'e bundle $P^{-1}$
under a composite morphism
$$\sU_{g-r}^r \stackrel {f_i} \longrightarrow \sU_{g-r} 
\stackrel \Delta \longrightarrow \sU_{g-r}\times_{\sA_{g-r}}\sU_{g-r},$$
where $\Delta$ is the diagonal embedding
and the rank $1$ projection $f_i$ depends upon $e_i$.
Since $\Delta^*P$ defines twice the principal polarization, \DHrefpart{i}
follows by taking $r=1$.

For \DHrefpart{ii}, we have now constructed $f_i:\sU_{g-r}^r\to \sU_{g-r}$
such that the class $H_i$ cuts out $f_i^*L_{g-r}$ on the generic fibre
of $\gamma:\sU_{g-r}^r\to \sA_{g-r}$. By \cite{FC}, p. 9,
$f_i$ extends uniquely to 
$\tpi_i:\sU_{g-r}^{r,part}\to \sU_{g-r}^{part}$.

%

Let $b,c\in\Q_{>0}$ and $n\in\N_{>0}$, to be determined later.
Put $\pi_i=[n]\circ\tpi_i$, where $[n]$ is multiplication
by $n$ on $\sU_{g-r}^{part}$.
We also let $\pi_i$ denote the induced morphism of geometric quotients
and $\rho_i=\alpha_{g-r+1,1}\circ\pi_i$.

Put $F=b\sum \rho_i^*(12M-D_{g-r+1}^{F})+c\gamma^*(12M-D_{g-r}^{part})$ and 
$G= \alpha_{g,r}^*(12M-D_{g})$.
It is enough, by \ref{H1},
to show that these classes are equal when restricted to the $0$-section 
$N_{g,r}^{part}$ of $\gamma$ and to the generic fibre $\Phi$ of $\gamma$.

On the zero-section
$N_{g,r}^{part}$, these restrictions are independent of $n$.
The restriction of $F$ is
$(pb+c)(12M-D_{g-r}^{part})$,
since $N_{g,r}^{part}$ maps isomorphically under $\pi_i$ to its image in
$U_{g-r}^{part}$, the $0$-section of $U_{g-r}^{part}\to A_{g-r}^{part}$.
On the other hand, $G\vert_{{N_{g,r}}^{part}}=12M-D_{g-r}^{part}$,
so the restrictions to $N_{g,r}^{part}$ are equal provided that
$c+p b=1$.

On $\Phi$, the class $M$ is trivial, so the restrictions are 
$$F\vert_\Phi = -bn^2\sum H_i,$$
since $[n]^*H_i = n^2H_i$, and
$$G\vert_\Phi = -s\sum H_i.$$
So we need to find $n,b,c$ such that $c+p b=1$ and $bn^2=s$.
Clearly this is possible.
\end{proof}
\end{proposition}

The next lemma removes the distinction between numerical and rational
equivalence in our context.

\begin{lemma} ${\tU}^r_{g-r}$ is regular. That is, its
Albanese variety is trivial.
\begin{proof} It is enough to prove this at some level $n\ge 3$ that is prime
to the characteristic. 

Put $U^{r,part}_{g-r,n} =U$ and $\tU^{r,part}_{g-r,n} =\tU$.
There is an open immersion $j:U\to \tU$ and $U$ is semi-abelian
over $A_{g-r,n}^{part}$. Then there is a commutative diagram
$$\begin{CD}
U @>>> A_{g-r,n}^{part}\cr
@VVV @VVV\cr
U^* @>{f}>> A_{g-r,n}^{Second}
\end{CD}$$
where the vertical arrows are open immersions, $f$ is proper and is
an equivariant compactification of a semi-abelian scheme over
$A_{g-r,n}^{Second}$. Over any point $P$ in $A_{g-r,n}^{Second}$
that lies over a $0$-dimensional cusp in $A_{g-r,n}^{Sat}$ the fibre
$f^{-1}(P)$ is stratified into (torsors under) algebraic tori, and so
is collapsed under the natural morphism
$\alpha:\Alb(U^*)\to \Alb(A_{g-r,n}^{Second})$. 
It follows that $\alpha$ is an isomorphism.

To prove the vanishing of $\Alb(A_{g-r,n}^{Second})$ it is enough,
by the usual comparison theorems, to work over $\C$. Suppose first
that $g-r\ge 2$. Then \cite{Kn} $A_{g-r,n}^{Second}$ is simply connected,
so regular. If $g-r=1$, then take $n=3$ or $5$ and use the fact that
$X(3)$ and $X(5)$ are rational.

Finally, the restriction homomorphism $\Pic \tU\to\Pic U$ is an isomorphism,
since $\tU\setminus U$ has codimension at least $2$, while
$\Pic U^*\to\Pic U$ is surjective.
\end{proof}
\end{lemma}

Now we can prove the main result of this section. This is just a
matter of assembling the pieces.

\begin{theorem}\label{Nef}
$12M-{D_g}$ is nef on $A_g^{F}$.

\begin{proof} If there is a curve $C$ in $A_g^{F}$ with $(12M-{D_g}).C<0$,
then for some $r\ge 2$ and some maximal cone $\sigma$ in ${\overline C}(X_r)$,
$C$ lies in the subvariety ${\overline U}_{g-r}^r$ corresponding to $\sigma$.
But \ref{modify} (replacing $g$ by $g-r$), \ref{product}
and \ref{codimension 2} together show that
the restriction of $12M-{D_g}$ to $\tU_{g-r}^r$ gives (after taking a
large multiple) a linear system whose base locus lies over $A_{g-r-1}^{Sat}$.
\end{proof}
\end{theorem}

\begin{corollary}
The cone $\NE(A_g^{F})$ of curves on $A_g^{F}$ is closed and
is the rational polyhedral cone generated by the curves
$C_1$ and $C_2$.

\begin{proof} 
Since $\rho(A_g^{F})=2$, 
by \ref{rho = 2}, it's enough to show that both curves lie in the 
boundary of $\cone(A_g^{F})$. 

For $C_2$ this is \ref{C_2 on boundary}.
For $C_1$ it follows from \ref{Nef} and the formula $(12M-{D_g}).C_1=0$.
\end{proof}
\end{corollary}

\begin{remark} As examples, consider the cases where $g=2$ and $n=2$ or $3$.

\noindent (i): $A^{F}_{2,2}$ is the blow-up of the Segre cubic threefold 
$S$ in its $10$ nodes.
The contraction $A_g^{F}=A^{F}_{2,2}/Sp_4(\F_2)\to S/Sp_4(\F_2)$ is the contraction
of the ray generated by $C_1$. 

\noindent (ii): $A^{F}_{2,3}$ is the blow-up of the Burkhardt quartic $B$ in 
its $45$ nodes and $A_g^{F}=A^{F}_{2,3}\to B/Sp_4(\F_3)$ is again
the contraction
of the ray generated by $C_1$. 
\end{remark}

\begin{corollary}\label{ample} \part[i]\label{HuS} 
The divisor class $aM-bD$ is ample on $A_g^{F}$ 
if and only if $12a>b>0$. 

\part[ii] $aM-bD$ is nef if and only if $12a\ge b\ge 0$.

\part[iii] The canonical class of $A_g^{F}$ is ample if and only if $g\ge 12$.
\begin{proof} Immediate from the result that $12M-D$ is nef
and Kleiman's criterion for ampleness:
a Cartier divisor class on the projective variety $X$ is ample if and 
only if it is strictly positive on $\partial\cone(X)$ (\cite{K}, Th. 2, p. 326).
\end{proof}
\end{corollary}

When $g\le 3$ \ref{HuS} is due to Hulek and Sankaran (\cite{HS02}, Theorem II.2.4).

\begin{definition}
The \emph{slope} of a Siegel cusp form is its weight divided by its order of vanishing.
(This is the reciprocal of what Weissauer \cite{W} calls the \emph{order of vanishing}
of a modular form.)
\end{definition}

\begin{corollary} Fix a slope $a\in\Q$. Then the ring of Siegel modular forms
of degree $g$ and slope $a$ 
and with Fourier coefficients in $\Z$ is a finitely generated
$\Z$-algebra provided that $a>12$. The same is true for $a=12$ provided that
$g\le 11$ and we consider Fourier coefficients in $\C$.

\begin{proof} For $a>12$ this follows from \ref{ample}. Over $\C$ it follows from
\ref{ample} and the base point free theorem for complex projective varieties 
with canonical singularities \cite{CKM}.
\end{proof}
\end{corollary}

Weissauer has shown (\cite{W}, p. 220) that for every $a>12$ and point $x\in A_g$,
there is a slope $a$ cusp form that does not vanish at $x$. This can be extended
to include $a=12$ (but not $a<12$, as he remarked).

\begin{corollary}\label{weissauer} For every $a \ge 12$ and every point $x\in A_g^{part}$
there is a slope $a$ cusp form that does not vanish (as a slope $a$ cusp form)
at $x$.

\begin{proof} When $a>12$ this follows from the ampleness on $A_g^{F}$
of the bundle of slope $a$ cusp forms. When $a=12$ we do not know, for
$g\ge 12$, whether
this bundle is eventually base-point-free on $A_g^{F}$; however, we do know,
by the proof of \ref{modify}, that it has no base points on $A_g^{part}$.
\end{proof}
\end{corollary}

\begin{remark} One natural problem is whether the total co-ordinate ring
of $A_g^{F}$, that is, the bigraded ring 
$\oplus_{a,b\ge 0} H^0(A_g^{F}, \sO(aM-bD))$
is of finite type, over $\Z$ or any other base ring. Another is to compute the 
various intersection numbers $M^a.D^b$ for $a+b=\frac{g(g+1)}{2}$.
\end{remark}

\end{section}
\bigskip

\begin{section}{Canonical models}\label{Tai}
\medskip
Over a field of characteristic zero
the results of the previous section
lead to statements about the field of Siegel modular functions. 
Over other fields the basic questions about resolving singularities,
even quotients of toroidal singularities by finite groups, are still too hard.
So in this section the base is $\Spec\C$. Identify $A_g$ with
$\frak H_g/\Gamma$, where $\frak H_g$ is the Siegel upper half-space
of degree $g$ and $\Gamma =Sp_{2g}(\Z)$. 
There is a birational
morphism $\pi:\tA\to\bA$.
(Since $\Gamma$ is not neat, $\tA$ will have non-trivial 
quotient singularities even if
the cone decompositions chosen for the construction of $\tA$
are basic. For quotients by neat arithmetic groups, these decompositions
are basic if and only if the toroidal compactification is smooth.) 

Let us say that a Deligne-Mumford stack $\sX$ has canonical or terminal singularities
if there is an \'etale surjective cover $X\to\sX$ from a scheme $X$ such that $X$
has canonical or terminal singularities.

\begin{lemma}\label{terminal stack} The stack $\sA_g^{F}$ has terminal singularities.

\begin{proof} Suppose that $T\inj X$ is an affine torus embedding.
This corresponds to a rational polyhedral convex cone $\sigma$ in 
$\mathbb X_*(T)\otimes \R$.
Then $X$ has $\Q$-Gorenstein singularities if and only if $\sigma$ is the cone
over a rational polyhedral convex polytope $\tau$ whose vertices are in
$\mathbb X_*(T)$ and that lies in an affine hyperplane $(z=1)$
of $\mathbb X_*(T)\otimes \R$. Moreover, $X$ has canonical singularities
if and only if in addition
there are no points in $\mathbb X_*(T)\cap\sigma$ satisfying $z<1$
and has terminal singularities if and only if the only points of 
$\mathbb X_*(T)\cap\sigma$
that satisfy $z\le 1$ are the vertices of $\tau$. 

To prove the result, we again use \ref{barnes}:
if $f$ is a quadratic form on $X_g$, then, when $f$ is
regarded as a linear function on ${\overline C}(X_g)$, 
its minima are all rank one
elements of $B(X_g)$. The statement now
follows from the description, which we have already recalled 
in Section \ref{voronoi}, 
of toroidal resolutions.
\end{proof}
\end{lemma}

To deal with the passage from the stack to its geometric quotient
requires analysis of the isotropy group actions. Tai shows \cite{T}
that the quotient map $\sA_g\to A_g$ is unramified in codimension $1$
when $g\ge 2$
and that $A_g$ has canonical singularities if $g\ge 5$; his argument
shows that in fact $A_g$ has terminal singularities if $g\ge 6$.
(Note that his proof of Lemma 4.3 of \emph{loc. cit.} needs slight amendment.
He states, without giving details, that
$$\sum\{\frac{t_i+t_j}{m}\}\ge \frac{r(r+1)^2}{4m},$$
where $\{x\}$ is the fractional part of $x$ and
$r=\half\phi(m)$; this leaves special cases, such as $m=30$, 
requiring individual treatment
beyond those that he gives. Attempts to fill in these details led to the estimate
$$\sum\{\frac{t_i+t_j}{m}\}\ge \frac{r^3+\frac{3}{2}r^2+\onehalf r}{4m};$$
this leaves the cases $m=8,10,12,18,30$ requiring individual treatment,
since then the last expression is less than $1$. Rather than check this by hand,
Tom Fisher wrote a Magma routine to do it. For $m=12$ he found
$\sum\{\frac{t_i+t_j}{m}\}=1$ and $\sum\{\frac{t_i+t_j}{m}\}> 1$ for the other values
of $m$.)

Tai's extension of this argument to deal with singularities in the boundary hints that
$A_g^{F}$ has canonical singularities if $g\ge 5$, and terminal singularities
if $g\ge 6$. Rather than explain this, we derive it from the facts that
$A_g$ and $\sA_g^F$ have canonical or terminal singularities and the following proposition,
which was proved by Snurnikov \cite{S} when $V$ is a point.

\begin{proposition}\label{snurnikov}
Suppose that $T\inj X$ is a torus embedding on which 
the finite group $G$ acts as a group of
algebraic torus automorphisms and that the $G$-action preserves
a certain neighbourhood $U$ of the boundary divisor $X\setminus T$. 
Put $U\cap T=U_0$.
Assume also that $V$ is a smooth variety on which $G$ acts,
that $G$ acts
freely in codimension $1$ on $U_0\times V$ and that $U$ and 
$(U_0\times V)/G$ have
canonical or terminal singularities. Then $G$ acts freely in codimension $1$
on $Z:=U\times V$ and $Z/G$ has canonical or terminal singularities.

\begin{proof} Choose a $T\rtimes G$-equivariant resolution $\tX\to X$. Let
$\tU$ denote the inverse image of $U$. We shall
show that $G$ acts freely in codimension $1$ on $\tZ:=\tU\times V$. If not, then 
there is a non-trivial subgroup $H$ of $G$ and an irreducible divisor $D$ on $\tZ$
such that $H$ acts trivially on $D$. Let $P\in D$ be general. Then
there is an $H$-equivariant isomorphism $T_{\tZ}(P)\to T_{\tU}(P)\oplus T_V(P)$.
So either $D=\tU\times D_2$ or $D=D_1\times V$. In the first case
$D\cap (U_0\times V)$ is non-empty, which is impossible, and so
$D=D_1\times V$, where $D_1$ is a boundary divisor. Then $D_1$ corresponds
to a 1-PS $\lambda:\GG_m\to T$ and the torus embedding $T\inj \tX$
gives a torus embedding $T_1:=T/\lambda(\GG_m)\inj D_1$ on which $H$ acts trivially.
But this contradicts the fact that $H$ acts freely in codimension $1$ on 
$U_0\times V$. So $G$ acts freely in codimension $1$ on both
$Z$ and $\tZ$.

Consider the commutative diagram
$$
\begin{CD}
\tZ @>>> \tZ/G\cr
@VVV @VVV\cr
Z @>>> Z/G.
\end{CD}
$$
For some integer $r\ge 1$ there is a $G$-invariant generator $\sigma$
of $\sO(rK_Z)$, which is then a generator of $\sO(rK_{Z/G})$.
Then $\sigma$ is either regular or zero along every exceptional divisor of
$\tZ\to Z$; since $G$ acts freely in codimension $1$ on $\tZ$,
it is enough to show that $\tZ/G$ has canonical or terminal singularities.
Then we can assume that $\tX=X=\A^n$ and that $T$ is the complement of the
co-ordinate hyperplanes, so that $G$ acts on $X$ by permuting the
co-ordinates. Let $L=\{(t,\ldots,t)\}$ be the diagonal copy of $\A^{F}$ in
$\A^n$. There is a $G$-equivariant decomposition $\A^n=L\times M$,
with $M=\A^{n-1}$, so that $(X\times V)/G$ is isomorphic to
$((M\times V)/G)\times L$. Since $(U_0\times V)/G$ has canonical or 
terminal singularities, we are done.
\end{proof}
\end{proposition}

\begin{corollary}\label{canonical} $A_g^{F}$ has canonical singularities if
$g\ge 5$ and terminal singularities if $g\ge 6$.

\begin{proof} According to the local description 
given immediately after \ref{decomposition},
over an arbitrary cusp, $\sA_g^{F}$ is a quotient stack
$[(X_1\times^{T_1}V)/GL(L_1)]$. Now the corollary is a consequence of 
\ref{snurnikov}, \ref{terminal stack}
and Tai's result, recalled above, that $A_g$ has canonical or terminal
singularities if $g\ge 5$ or $g\ge 6$.
\end{proof}
\end{corollary}

\begin{corollary} \part[i] $A_g^{F}\to A_g^{Sat}$ is the relative canonical model
of $A_g^{Sat}$ if $g\ge 5$.

\part[ii] $A_g^{F}$ is the canonical model of $A_g$ if $g\ge 12$.

\part[iii] The canonical model of $A_{11}$ exists and
is the result of contracting the extremal ray $\R_+[C_1]$.

\begin{proof} \DHrefpart{i} and \DHrefpart{ii} 
follow from \ref{ample} and \ref{canonical}.
For \DHrefpart{iii}, we also need the base point free theorem.
Note that since $K_{A_{11}}.C_2 > 0$, the result of the contraction, which is guaranteed
to have terminal singularities, has ample canonical class.
\end{proof}
\end{corollary}

In particular, we recover the result, weaker than that of Freitag, Mumford and Tai,
that $A_g$ is of general type if $g\ge 11$. If $5\le g\le 10$, then $C_1$
spans an extremal ray on which $K$ is negative. This ray is then contractible,
from general results on complex varieties. It would be interesting to understand
how to carry forward the minimal model program for these varieties, especially
when $g=6$, since this is the only case where the Kodaira dimension of $A_g$
is still unknown.

\end{section}
\bigskip
\begin{section}{Higher level}\label{higher}
\medskip
In this section we work over a field $k$ that contains $\zeta_n$, 
and $n$ is an integer with $n\ge 2$ that is not divisible by $\ch k$.

\begin{theorem} 
\part[i] $aM-D^{(n)}$ is nef on the first Voronoi compactification $A^{F}_{g,n}$ of
$A_{g,n}$, where $D^{(n)}$ is the reduced boundary divisor, if and only if $a\ge 12/n$
and ample if and only if $a> 12/n$. 

\part[ii] The canonical class of $A^{F}_{g,n}$ is ample if $g+1> 12/n$.

\part[iii] Assume that $\ch k=0$ and that either $g\ge 3$ or $n\ge 3$. 
Then $A^{F}_{g,n}$ has canonical singularities
and is the relative canonical model of $A_{g,n}^{Sat}$. Moreover,
$A_{g,n}$ is of general type if $g+1\ge 12/n$ and $A_{g,n}^F$ is its 
canonical model if $g+1> 12/n$.

\begin{proof} Consider the projection $\pi:A^{F}_{g,n}\to A_g^{F}$. Then
$\pi^*D=nD^{(n)}$, so that $K_{A^{F}_{g,n}}=\pi^*((g+1)M-\frac{1}{n}D)$, 
and now \DHrefpart{i} and \DHrefpart{ii} follow
from the description of the ample cone on $A_g^{F}$.
The proof of \DHrefpart{iii} is also immediate.
\end{proof}
\end{theorem}

\begin{remark} In fact, $A_{g,n}$ is known to be of general type in this range,
and more besides; see \cite{HS02}, Theorem II.2.1, p. 106.

\end{remark}
\end{section}
\bigskip
\bibliography{alggeom,ekedahl}

\providecommand{\bysame}{\leavevmode\hbox to3em{\hrulefill}\thinspace}
\begin{thebibliography}{EGAIII:2}

\bibitem[A]{A}
V.~Alexeev, \emph{Complete moduli in the presence of semi-abelian group action},
Annals of Math. \textbf{155} (2002), 611-708.

\bibitem[AMRT]{AMRT} A.~Ash, D.~Mumford, M.~Rapoport and Y.-S.~Tai,
\emph{Smooth compactifications of locally symmetric varieties}, Math Sci Press,
Boston, 1975.

\bibitem[BC]{BC} E.S.~Barnes and M.J.~Cohn, \emph{On the inner product of positive
quadratic forms}, J. London Math. Soc. \textbf{12} (1975/76), 32-36.

\bibitem[BMS]{BMS}
H.~Bass, J.~Milnor and J.-P.~Serre, \emph{Solution of the congruence subgroup
problem}, Pub. Math. IHES \textbf{33} (1967), 59-137.

\bibitem[CKM]{CKM}
H.~Clemens, J.~Koll\'ar and S.~Mori, \emph{Higher dimensional complex geometry},
Ast\'erisque \textbf{166}, 1988.

\bibitem[ER]{ER}
R.~Erdahl and K.~Rybnikov, Voronoi--Dickson hypothesis on perfect forms and $L$-types,
math.NT/0112097.

\bibitem[FC]{FC}
G.~Faltings and C.~L. Chai, \emph{Degenerations of abelian varieties},
Springer Verlag, 1980.

\bibitem[F1]{F1}
E.~Freitag, \emph{Siegelsche Modulfunktionen}, Springer Verlag, 1983.

\bibitem[F2]{F2}
E.~Freitag, \emph{Die Irreduzibilit\"at der Schottkyrelation (Bemerkung \"uber
ein Satz von J.~Igusa)}, Arch. Math. \textbf{40} (1983), 255-259.

\bibitem[vdG]{vdG}
G.~van der Geer, \emph{Hilbert modular surfaces}, Springer Verlag, 1988.

\bibitem[HS02]{HS02}
K.~Hulek and G.K.~Sankaran, \emph{The geometry of Siegel modular varieties},
in Higher Dimensional Birational Geometry, Adv. Studies in Pure Math. \textbf{35},
Tokyo, 2002.

\bibitem[HS04]{HS04}
K.~Hulek and G.K.~Sankaran, \emph{The nef cone of toroidal compactifications of $A_4$},
J. London Math. Soc. , \textbf{88} (2004), 659-704.

\bibitem [K]{K}
S.~Kleiman, \emph{Towards a numerical theory of ampleness}, Annals of Math. \textbf{84}
(1966), 293-344.

\bibitem[KM]{KM}
S.~Keel and S.~Mori, \emph{Quotients by groupoids}, Annals of Math. \textbf{145}
(1997), 193-213.

\bibitem [Kn]{Kn}
F.W.~Kn\"oller, \emph{Die Fundamentalgruppen der Siegelsche Modulvariet\"aten},
Abh. Math. Seminar Univ. Hamburg \textbf{57} (1987), 203-213.

\bibitem[LB]{LB} H.~Lange and C.~Birkenhake, \emph{Complex abelian varieties},
Springer, Berlin, 1992.

%
\bibitem[M]{M}
D.~Mumford, \emph{On the Kodaira dimension of the Siegel modular variety},
Lecture Notes
in Mathematics \textbf{997}, Springer, Berlin, 1983.

\bibitem[N]{N}
Y.~Namikawa, \emph{Toroidal compactification of Siegel spaces}, Lecture Notes
in Mathematics \textbf{812}, Springer, Berlin, 1980.

\bibitem[S]{S} V.S.~Snurnikov, Ph.D. thesis, Cambridge, 2002.

\bibitem[T]{T} Y.-S.~Tai, \emph{On the Kodaira dimension of the moduli space 
of abelian varieties}, Invent. Math. \textbf{68} (1982), 425-439.

\bibitem[W]{W} R.~Weissauer, \emph{Untervariet\"aten der Siegelschen Modulmannigfaltigkeiten
von allgemeinem Typ}, Math. Annalen \textbf{275} (1986), 207-220.

\end{thebibliography}
\bibliographystyle{pretex}
\end{document}